\documentclass[11pt]{amsart}
\usepackage{mathrsfs}
\usepackage{amssymb,latexsym}
\usepackage{pstcol,pstricks,color}
\usepackage{chemarrow}

 \numberwithin{equation}{section}

\newtheorem{theorem}{Theorem}[section]
\newtheorem{proposition}[theorem]{Proposition}

\newtheorem{lemma}[theorem]{Lemma}

\def\qed{\hfill $\Box$}
\def\pf{\noindent {\it Proof.} }

\title[]{A bijection between the sets of $(a,b,b^2)$-Generalized Motzkin paths avoiding $\mathbf{uvv}$-patterns and $\mathbf{uvu}$-patterns }

\begin{document}
\maketitle

\begin{center}
Yidong Sun\footnote{Corresponding author: Yidong Sun.}, Cheng Sun$^{2}$ and Xiuli Hao$^{3}$

School of Science, Dalian Maritime University, 116026 Dalian, P.R. China\\[5pt]

{\it Emails: $^{1}$sydmath@dlmu.edu.cn, $^{2}$scmath@dlmu.edu.cn, $^{3}$hxl1212032018@163.com}

\end{center}\vskip0.2cm

\subsection*{Abstract} A generalized Motzkin path, called G-Motzkin path for short, of length $n$ is a lattice path from $(0, 0)$ to $(n, 0)$ in the first quadrant of the XOY-plane that consists of up steps $\mathbf{u}=(1, 1)$, down steps $\mathbf{d}=(1, -1)$, horizontal steps $\mathbf{h}=(1, 0)$ and vertical steps $\mathbf{v}=(0, -1)$. An $(a,b,c)$-G-Motzkin path is a weighted G-Motzkin path such that the $\mathbf{u}$-steps, $\mathbf{h}$-steps, $\mathbf{v}$-steps and $\mathbf{d}$-steps are weighted respectively by $1, a, b$ and $c$.
Let $\tau$ be a word on $\{\mathbf{u}, \mathbf{d}, \mathbf{v}, \mathbf{d}\}$, denoted by $\mathcal{G}_n^{\tau}(a,b,c)$ the set of $\tau$-avoiding $(a,b,c)$-G-Motzkin paths of length $n$ for a pattern $\tau$. In this paper, we consider the $\mathbf{uvv}$-avoiding $(a,b,c)$-G-Motzkin paths and provide a direct bijection $\sigma$ between $\mathcal{G}_n^{\mathbf{uvv}}(a,b,b^2)$ and $\mathcal{G}_n^{\mathbf{uvu}}(a,b,b^2)$. Finally, the set of fixed points of $\sigma$ is also described and counted.

\medskip

{\bf Keywords}: G-Motzkin path, Catalan number.

\noindent {\sc 2020 Mathematics Subject Classification}: Primary 05A15, 05A19; Secondary 05A10.

{\bf \section{ Introduction } }

A {\it generalized Motzkin path}, called {\it G-Motzkin path} for short, of length $n$ is a lattice path from $(0, 0)$ to $(n, 0)$ in the first quadrant of the XOY-plane that consists of up steps $\mathbf{u}=(1, 1)$, down steps $\mathbf{d}=(1, -1)$, horizontal steps $\mathbf{h}=(1, 0)$ and vertical steps $\mathbf{v}=(0, -1)$.  Other related lattice paths with various steps including vertical steps permitted have been considered by \cite{Dziem-A, Dziem-B, Dziem-C, YanZhang, SunZhao, SunWang}. See Figure 1 for a G-Motzkin path of length $25$.

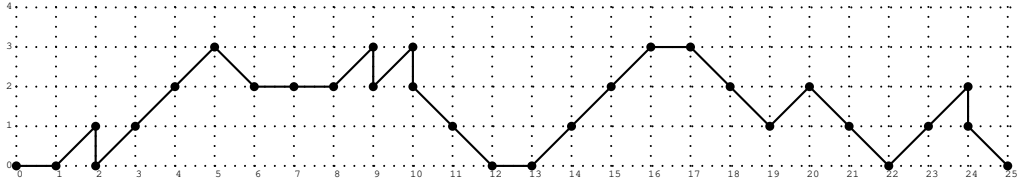
\begin{figure}[h] \setlength{\unitlength}{0.5mm}

\begin{center}
\begin{pspicture}(13,2.2)
\psset{xunit=15pt,yunit=15pt}\psgrid[subgriddiv=1,griddots=4,
gridlabels=4pt](0,0)(25,4)

\psline(0,0)(1,0)(2,1)(2,0)(3,1)(5,3)(6,2)(8,2)(9,3)(9,2)(10,3)(10,2)(11,1)(12,0)(13,0)(16,3)(17,3)(18,2)(19,1)
\psline(19,1)(20,2)(22,0)(24,2)(24,1)(25,0)

\pscircle*(0,0){0.06}\pscircle*(1,0){0.06}\pscircle*(2,1){0.06}\pscircle*(2,0){0.06}
\pscircle*(3,1){0.06}\pscircle*(4,2){0.06}\pscircle*(5,3){0.06}\pscircle*(6,2){0.06}

\pscircle*(7,2){0.06}\pscircle*(8,2){0.06}\pscircle*(9,2){0.06}
\pscircle*(9,3){0.06}\pscircle*(10,2){0.06}\pscircle*(10,3){0.06}\pscircle*(11,1){0.06}
\pscircle*(12,0){0.06}\pscircle*(13,0){0.06}\pscircle*(14,1){0.06}
\pscircle*(15,2){0.06}\pscircle*(16,3){0.06}\pscircle*(17,3){0.06}\pscircle*(18,2){0.06}
\pscircle*(19,1){0.06}\pscircle*(20,2){0.06}
\pscircle*(21,1){0.06}\pscircle*(22,0){0.06}\pscircle*(23,1){0.06}\pscircle*(24,2){0.06}
\pscircle*(24,1){0.06}\pscircle*(25,0){0.06}

\end{pspicture}
\end{center}\vskip0.2cm

\caption{\small A G-Motzkin path of length $25$.}

\end{figure}

An {\it $(a,b,c)$-G-Motzkin path} is a weighted G-Motzkin path $\mathbf{P}$ such that the $\mathbf{u}$-steps, $\mathbf{h}$-steps, $\mathbf{v}$-steps and $\mathbf{d}$-steps of $\mathbf{P}$ are weighted respectively by $1, a, b$ and $c$. The {\it weight} of $\mathbf{P}$, denoted by $w(\mathbf{P})$, is the product of the weight of each step of $\mathbf{P}$. For example, $w(\mathbf{uhuduuvvdhh})=a^3b^2c^2$. The {\em weight} of a subset $\mathcal{A}$ of the set of weighted G-Motzkin paths, denoted by $w(\mathcal{A})$, is the sum of the total weights of all paths in $\mathcal{A}$. Denoted by $G_n(a, b, c)$ the weight of the set $\mathcal{G}_{n}(a, b, c)$ of all $(a,b,c)$-G-Motzkin paths of length $n$.
Let $\tau$ be a word on $\{\mathbf{u}, \mathbf{d}, \mathbf{v}, \mathbf{d}\}$, denoted by $G_n^{\tau}(a, b, c)$ the weight of the set $\mathcal{G}_{n}^{\tau}(a, b, c)$ of all $\tau$-avoiding $(a,b,c)$-G-Motzkin paths of length $n$, that is the weight of the subset of all $(a,b,c)$-G-Motzkin paths of length $n$ avoiding the pattern $\tau$. Figure 1 is an example of a G-Motzkin paths of length $25$ avoiding the pattern $\mathbf{uvv}$, but not avoiding the pattern $\mathbf{uvu}$.

Recently Sun et al. \cite{SunZhao, SunWang} have derived the generating functions of $G_n(a, b, c)$ and $G_n^{\mathbf{uvu}}(a, b, c)$ as follows
\begin{eqnarray}
G(a,b,c; x) \hskip-.22cm &=&\hskip-.22cm \sum_{n=0}^{\infty}G_n(a, b, c)x^n     \nonumber \\
            \hskip-.22cm &=&\hskip-.22cm \frac{1-ax-\sqrt{(1-ax)^2-4x(b+cx)}}{2x(b+cx)}=\frac{1}{1-ax}C\Big(\frac{x(b+cx)}{(1-ax)^2}\Big),  \label{eqn 1.1}\\
G^{\mathbf{uvu}}(a,b,c; x) \hskip-.22cm &=&\hskip-.22cm \sum_{n=0}^{\infty}G_n^{\mathbf{uvu}}(a, b, c)x^n \nonumber \\
                           \hskip-.22cm &=&\hskip-.22cm \frac{(1-ax)(1+bx)-\sqrt{(1-ax)^2(1+bx)^2-4x(1+bx)(b+cx)}}{2x(b+cx)} \label{eqn 1.2} \\
                           \hskip-.22cm &=&\hskip-.22cm  \frac{1}{1-ax}C\Big(\frac{x(b+cx)}{(1-ax)^2(1+bx)}\Big), \nonumber
\end{eqnarray}
where
\begin{eqnarray}\label{eqn 1.3}
C(x)=\sum_{n=0}^{\infty}C_nx^n=\frac{1-\sqrt{1-4x}}{2x}
\end{eqnarray}
is the generating function for the well-known Catalan number $C_n=\frac{1}{n+1}\binom{2n}{n}$, counting the number of Dyck paths of length of $2n$ \cite{Stanley, StanleyEC2}.

A {\it Dyck path} of length $2n$ is a G-Motzkin path of length $2n$ with no $\mathbf{h}$-steps or $\mathbf{v}$-steps. A {\it Motzkin path} of length $n$ is a G-Motzkin path of length $n$ with no $\mathbf{v}$-steps. An {\it $(a,b)$-Dyck path} is a weighted Dyck path with $\mathbf{u}$-steps weighted by $1$, $\mathbf{d}$-steps in $\mathbf{ud}$-peaks weighted by $a$ and other $\mathbf{d}$-steps weighted by $b$. An {\it $(a,b)$-Motzkin path} of length $n$ is an $(a,0,b)$-G-Motzkin path of length $n$.
A {\it Schr\"{o}der path} of length $2n$ is a path from $(0, 0)$ to $(2n, 0)$ in the first quadrant of the XOY-plane that consists of up steps $\mathbf{u}=(1, 1)$, down steps $\mathbf{d}=(1, -1)$ and horizontal steps $\mathbf{H}=(2, 0)$. An {\it $(a,b)$-Schr\"{o}der path} is a weighted Schr\"{o}der path such that the $\mathbf{u}$-steps, $\mathbf{H}$-steps and $\mathbf{d}$-steps are weighted respectively by $1, a$ and $b$.

Let $\mathcal{C}_n(a,b), \mathcal{M}_n(a,b)$ and $\mathcal{S}_n(a,b)$ be respectively the sets of $(a,b)$-Dyck paths of length $2n$, $(a,b)$-Motzkin paths of length $n$ and $(a,b)$-Schr\"{o}der paths of length $2n$. Let $C_n(a,b), M_n(a,b)$ and $S_n(a,b)$ be their weights with $C_0(a,b)=M_0(a,b)=S_0(a,b)=1$ respectively.
It is not difficult to deduce that \cite{ChenPan}
\begin{eqnarray*}
C_n(a,b)  \hskip-.22cm &=&\hskip-.22cm  \sum_{k=1}^{n}\frac{1}{n}\binom{n}{k-1}\binom{n}{k}a^{k}b^{n-k},  \\
M_n(a,b)  \hskip-.22cm &=&\hskip-.22cm  \sum_{k=0}^{n}\binom{n}{2k}C_ka^{n-2k}b^{k},                  \\
S_n(a,b)  \hskip-.22cm &=&\hskip-.22cm  \sum_{k=0}^{n}\binom{n+k}{2k}C_ka^{n-k}b^{k},
\end{eqnarray*}
and their generating functions
\begin{eqnarray*}
C(a,b; x)\hskip-.22cm &=&\hskip-.22cm \sum_{n=0}^{\infty}C_n(a,b)x^n=\frac{1-(a-b)x-\sqrt{(1-(a-b)x)^2-4bx}}{2bx}, \\
M(a,b; x)\hskip-.22cm &=&\hskip-.22cm \sum_{n=0}^{\infty}M_n(a,b)x^n=\frac{1-ax-\sqrt{(1-ax)^2-4bx^2}}{2bx^2}, \\
S(a,b; x)\hskip-.22cm &=&\hskip-.22cm \sum_{n=0}^{\infty}S_n(a,b)x^n=\frac{1-ax-\sqrt{(1-ax)^2-4bx}}{2bx}.
\end{eqnarray*}
There are closely relations between $C_n(a,b), M_n(a,b)$ and $S_n(a,b)$. Exactly, Chen and Pan \cite{ChenPan} derived the following equivalent relations
\begin{eqnarray*}
S_n(a,b)=C_n(a+b,b)=(a+b)M_{n-1}(a+2b,(a+b)b)
\end{eqnarray*}
for $n\geq 1$ and provided some combinatorial proofs. Sun et al. \cite{SunWang} obtained that
\begin{eqnarray*}
G_{n}^{\mathbf{uvu}}(a, b, b^2)=S_n(a,b)
\end{eqnarray*}
for $n\geq 0$ and presented bijections between the sets $\mathcal{G}_{n}^{\mathbf{uvu}}(a, b, b^2)$  and $\mathcal{S}_n(a,b)$ as well as the set $\mathcal{C}_n(a+b,b)$.

In the literature, there are many papers dealing with $(a,b)$-Motzkin paths. For examples, Chen and Wang \cite{ChenWang} explored the connection between noncrossing linked partitions and $(3,2)$-Motzkin paths, established a one-to-one correspondence between the set of noncrossing linked partitions of $\{1, \dots, n+1\}$ and the set of large $(3,2)$-Motzkin paths of length $n$, which leads to a simple explanation of the well-known relation between the large and the little Schr\"{o}der numbers. Yan \cite{Yan07} found a bijective proof between the set of restricted $(3,2)$-Motzkin paths of length $n$ and the set of the Schr\"{o}der paths of length $2n$.

In the present paper we concentrate on the $\mathbf{uvv}$-avoiding G-Motzkin paths, that is, the G-Motzkin paths with no $\mathbf{uvv}$ patterns. Precisely, the next section considers the enumerations of the set of $\mathbf{uvv}$-avoiding $(a,b,c)$-G-Motzkin paths and the set of $\mathbf{uvv}$-avoiding $(a,b,c)$-G-Motzkin paths with no $\mathbf{h}$-steps on the $x$-axis, and find that $G_{n}^{\mathbf{uvv}}(a, b, b^2)=G_{n}^{\mathbf{uvu}}(a, b, b^2)$. The third section provides a direct bijection $\sigma$ between the set $\mathcal{G}_n^{\mathbf{uvv}}(a,b,b^2)$ of $\mathbf{uvv}$-avoiding $(a,b,b^2)$-G-Motzkin paths and the set $\mathcal{G}_n^{\mathbf{uvu}}(a,b,b^2)$ of $\mathbf{uvu}$-avoiding $(a,b,b^2)$-G-Motzkin paths. Finally, the set of fixed points of $\sigma$ is also described and counted.

\section{ $\mathbf{uvv}$-avoiding $(a,b,c)$-G-Motzkin paths }

In this section, we first consider the $\mathbf{uvv}$-avoiding $(a,b,c)$-G-Motzkin paths which involve some classical structures as special cases, and count the set of $\mathbf{uvv}$-avoiding $(a,b,c)$-G-Motzkin paths with no $\mathbf{h}$-steps on the $x$-axis.

Let $G^{\mathbf{uvv}}(a,b,c; x)=\sum_{n=0}^{\infty}G_n^{\mathbf{uvv}}(a, b, c)x^n$ be the generating function for the $\mathbf{uvv}$-avoiding $(a,b,c)$-G-Motzkin paths. According to the method of the first return decomposition \cite{Deutsch99}, any $\mathbf{uvv}$-avoiding $(a,b,c)$-G-Motzkin path $\mathbf{P}$ can be decomposed as one of the following four forms:
$$\mathbf{P}=\varepsilon, \ \mathbf{P}=\mathbf{h}_a\mathbf{Q}_1, \ \mathbf{P}=\mathbf{u}\mathbf{Q}_2\mathbf{d}_c\mathbf{Q}_1 \ \mbox{or}\  \mathbf{P}=\mathbf{u}\mathbf{Q}_3\mathbf{v}_b\mathbf{Q}_1, $$
where $\mathbf{x}_t$ denotes the $\mathbf{x}$-steps with weight $t$, $\mathbf{Q}_1$ and $\mathbf{Q}_2$ are (possibly empty) $\mathbf{uvv}$-avoiding $(a,b,c)$-G-Motzkin paths, and $\mathbf{Q}_3$ is any $\mathbf{uvv}$-avoiding $(a,b,c)$-G-Motzkin paths with no $\mathbf{uv}$-step at the end of $\mathbf{Q}_3$. Then we get the relation
\begin{eqnarray}
G^{\mathbf{uvv}}(a,b,c; x) \hskip-.22cm &=&\hskip-.22cm  1+axG^{\mathbf{uvv}}(a,b,c; x)+cx^2G^{\mathbf{uvv}}(a,b,c; x)^2         \nonumber \\
                           \hskip-.22cm & &\hskip-.22cm  \ \ +\ bx\big(G^{\mathbf{uvv}}(a,b,c; x)-bxG^{\mathbf{uvv}}(a,b,c; x)\big)G^{\mathbf{uvv}}(a,b,c; x)      \nonumber\\
                           \hskip-.22cm &=&\hskip-.22cm  1+axG^{\mathbf{uvv}}(a,b,c; x)+(b+(c-b^2)x)xG^{\mathbf{uvv}}(a,b,c; x)^2.                           \label{eqn 2.1}
\end{eqnarray}
Solve this, we have
\begin{eqnarray}\label{eqn 2.2}
G^{\mathbf{uvv}}(a,b,c; x) \hskip-.22cm &=&\hskip-.22cm \frac{1-ax-\sqrt{(1-ax)^2-4x(b+(c-b^2)x)}}{2x(b+(c-b^2)x)} \nonumber \\
                          \hskip-.22cm &=&\hskip-.22cm  \frac{1}{1-ax}C\Big(\frac{x(b+(c-b^2)x)}{(1-ax)^2}\Big).
\end{eqnarray}
When $a=b=c=1$, $G^{\mathbf{uvv}}(1,1,1; x)=\frac{1-x-\sqrt{(1-x)^2-4x}}{2x}$ is just the generating function of the large Schr\"{o}der numbers \cite{Sloane}.

By (\ref{eqn 1.3}), taking the coefficient of $x^n$ in $G^{\mathbf{uvv}}(a,b,c; x)$, we derive that
\begin{proposition}
For any integer $n\geq 0$, there holds
\begin{eqnarray*}
G_n^{\mathbf{uvv}}(a,b,c) \hskip-.22cm &=&\hskip-.22cm \sum_{k=0}^{n}\sum_{j=0}^{k}\sum_{\ell=0}^{n-k-j}(-1)^{\ell}\binom{k}{j}\binom{k+\ell-1}{\ell}\binom{n+k-j-\ell}{2k}C_ka^{n-k-j-\ell}b^{k+\ell-j}(c-b^2)^{j}  \\
                          \hskip-.22cm &=&\hskip-.22cm \sum_{k=0}^{n}\sum_{j=0}^{k}\sum_{\ell=0}^{n-k-j}(-1)^{n-k-j-\ell}\binom{k}{j}\binom{2k+\ell}{\ell}\binom{n-j-\ell-1}{n-k-j-\ell}C_ka^{\ell}b^{n-2j-\ell}(c-b^2)^{j}.
\end{eqnarray*}
\end{proposition}

Set $T=xG(a,b,c; x)$, (\ref{eqn 2.1}) produces
\begin{eqnarray}\label{eqn 2.3}
T=x\frac{1+aT+(c-b^2)T^2}{1-bT},
\end{eqnarray}
using the Lagrange inversion formula \cite{Gessel}, taking the coefficient of $x^{n+1}$ in $T$ in three different ways, we derive that
\begin{proposition}
For any integer $n\geq 0$, there holds
\begin{eqnarray*}
G_n^{\mathbf{uvv}}(a,b,c) \hskip-.22cm &=&\hskip-.22cm \frac{1}{n+1}\sum_{k=0}^{[\frac{n}{2}]}\sum_{j=0}^{n-2k}\binom{n+1}{k}\binom{n+1-k}{j}\binom{2n-2k-j}{n-2k-j} a^{j}b^{n-2k-j} (c-b^2)^{k}  \nonumber\\
             \hskip-.22cm &=&\hskip-.22cm \frac{1}{n+1}\sum_{k=0}^{n}\sum_{j=0}^{[\frac{n-k}{2}]}\binom{n+1}{k}\binom{n+1-k}{j}\binom{2n-k-2j}{n-k-2j} a^{k}b^{n-k-2j}(c-b^2)^{j}  \nonumber\\
             \hskip-.22cm &=&\hskip-.22cm \frac{1}{n+1}\sum_{k=0}^{n}\sum_{j=0}^{n-k}\binom{n+1}{k}\binom{k}{j}\binom{2n-k-j}{n-k-j} a^{k-j}b^{n-k-j}(c-b^2)^{j}.
\end{eqnarray*}
\end{proposition}

Exactly, by (\ref{eqn 1.1}), (\ref{eqn 1.2}) and (\ref{eqn 2.2}), it can be deduced that
$$ G^{\mathbf{uvv}}(a,b,b^2+c; x)= G(a,b,c; x), \ \ G^{\mathbf{uvv}}(a,b,b^2; x)= G^{\mathbf{uvu}}(a,b,b^2; x). $$
That is $G_n^{\mathbf{uvv}}(a,b,b^2+c)=G_n(a,b,c)$ and $G_n^{\mathbf{uvv}}(a,b,b^2)=G_n^{\mathbf{uvu}}(a,b,b^2)$. The first identity has a direct combinatorial interpretation if one notices that each $\mathbf{d}_{b^2+c}$-step of $\mathbf{P}\in \mathcal{G}_{n}^{\mathbf{uvv}}(a, b, b^2+c)$ can be regarded equivalently as the corresponding $\mathbf{d}_{c}$-step and $\mathbf{u}\mathbf{v}_{b}\mathbf{v}_{b}$-step of $\mathbf{P}'\in \mathcal{G}_{n}(a, b, c)$ . The combinatorial interpretation of the second identity will be given in the next section.

When $(a,b,c)$ is specialized, $G^{\mathbf{uvv}}(a,b,c; x)$ and $G_n^{\mathbf{uvv}}(a,b,c)$ reduce to some well-known generating functions and classical combinatorial sequences involving the Catalan numbers $C_n$, Motzkin numbers $M_n$, the large Schr\"{o}der numbers $S_n$, $(a+b,b)$-Catalan number $C_n(a+b,b)$, $(a,b)$-Motzkin number $M_n(a,b)$ and $(a,b)$-Schr\"{o}der number $S_n(a,b)$. See Table 2.1 for example.
\begin{center}
\begin{eqnarray*}
\begin{array}{c|c|c|c}\hline
 (a, b, c)   & G^{\mathbf{uvv}}(a,b,c; x)               &     G_n^{\mathbf{uvv}}(a,b,c)      &       Senquences                         \\[5pt]\hline
 (0, 1, 1)   & C(x)=\frac{1-\sqrt{1-4x}}{2x}            &     C_n                            &       \mbox{\cite [A000108]{Sloane} }    \\[5pt]\hline
 (1, 0, 1)   & M(x)=\frac{1-x-\sqrt{1-2x-3x^2}}{2x^2}   &     M_n                            &       \mbox{\cite [A001006]{Sloane} }     \\[5pt]\hline
 (1, 1, 1)   & S(x)=\frac{1-x-\sqrt{1-6x+x^2}}{2x}      &     S_n                            &       \mbox{\cite [A006318]{Sloane} }      \\[5pt]\hline
 (1, 0, 2)   & \frac{1-x-\sqrt{1-2x-7x^2}}{4x}          &     a_n                            &       \mbox{\cite [A025235]{Sloane} }       \\[5pt]\hline
 (-3, 4, 16) & \frac{1+3x-\sqrt{1-10x+9x^2}}{8x}        &     a_n                            &       \mbox{\cite [A059231]{Sloane} }      \\[5pt]\hline
 (a, 0, b)   & \frac{1-ax-\sqrt{(1-ax)^2-4bx^2}}{2bx^2} &    M_n(a,b)                        &          \\[5pt]\hline
 (a, b, b^2) & \frac{1-ax-\sqrt{(1-ax)^2-4bx}}{2bx}     &    C_n(a+b,b)\ or\ S_n(a,b)        &          \\[5pt]\hline
\end{array}
\end{eqnarray*}
Table 2.1. The specializations of $G^{\mathbf{uvv}}(a,b,c; x)$ and $G_n^{\mathbf{uvv}}(a,b,c)$.
\end{center}

Denoted by $\bar{G}_n^{\mathbf{uvv}}(a, b, c)$ the weight of the set $\mathcal{\bar{G}}_{n}^{\mathbf{uvv}}(a, b, c)$ of all $\mathbf{uvv}$-avoiding $(a,b,c)$-G-Motzkin paths of length $n$ such that the paths have no $\mathbf{h}$-steps on the $x$-axis. Set $\mathcal{\bar{G}}^{\mathbf{uvv}}(a, b, c)=\bigcup_{n\geq 0}\mathcal{\bar{G}}_n^{\mathbf{uvv}}(a, b, c)$.

Let $\bar{G}^{\mathbf{uvv}}(a,b,c; x)=\sum_{n=0}^{\infty}\bar{G}_n^{\mathbf{uvv}}(a, b, c)x^n$ be the generating function for the $\mathbf{uvv}$-avoiding $(a,b,c)$-G-Motzkin paths in $\mathcal{\bar{G}}^{\mathbf{uvv}}(a, b, c)$. According to the method of the first return decomposition, any paths $\mathbf{P}\in \mathcal{\bar{G}}^{\mathbf{uvv}}(a, b, c)$ can be decomposed as one of the following three forms:
$$\mathbf{P}=\varepsilon, \ \mathbf{P}=\mathbf{u}\mathbf{Q}_2\mathbf{d}_c\mathbf{Q}_1 \ \mbox{or}\  \mathbf{P}=\mathbf{u}\mathbf{Q}_3\mathbf{v}_b\mathbf{Q}_1, $$
where $\mathbf{Q}_1\in \mathcal{\bar{G}}^{\mathbf{uvv}}(a, b, c)$, $\mathbf{Q}_2\in \mathcal{G}^{\mathbf{uvv}}(a, b, c)$ and $\mathbf{Q}_3\in \mathcal{{G}}^{\mathbf{uvv}}(a, b, c)$ has no $\mathbf{uv}$-step at the end of $\mathbf{Q}_3$. Then we get the relation
\begin{eqnarray*}
\bar{G}^{\mathbf{uvv}}(a,b,c; x) \hskip-.22cm &=&\hskip-.22cm  1+cx^2G^{\mathbf{uvv}}(a,b,c; x)\bar{G}^{\mathbf{uvv}}(a,b,c; x)         \nonumber \\
                           \hskip-.22cm & &\hskip-.22cm  \ \ +\ bx\big(G^{\mathbf{uvv}}(a,b,c; x)-bxG^{\mathbf{uvv}}(a,b,c; x)\big)\bar{G}^{\mathbf{uvv}}(a,b,c; x),
\end{eqnarray*}
which, by (\ref{eqn 2.1}) and (\ref{eqn 2.3}), leads to
\begin{eqnarray*}
x\bar{G}^{\mathbf{uvv}}(a,b,c; x) \hskip-.22cm &=&\hskip-.22cm  \frac{x}{1-(b+(c-b^2)x)xG^{\mathbf{uvv}}(a,b,c; x)}=\frac{xG^{\mathbf{uvv}}(a,b,c; x)}{1+axG^{\mathbf{uvv}}(a,b,c; x)}=\frac{T}{1+aT}.
\end{eqnarray*}
By the Lagrange inversion formula, taking the coefficient of $x^{n+1}$ in $x\bar{G}^{\mathbf{uvv}}(a,b,c; x)$ in three different ways, we derive that
\begin{proposition}
For any integer $n\geq 0$, there holds
\begin{eqnarray*}
\bar{G}_n^{\mathbf{uvv}}(a,b,c) \hskip-.22cm &=&\hskip-.22cm \sum_{i=0}^{n+1}(-1)^i\frac{i+1}{n+1}\sum_{k=0}^{[\frac{n}{2}]}\sum_{j=0}^{n-2k}\binom{n+1}{k}\binom{n+1-k}{j} \\
                                \hskip-.22cm & &\hskip-.22cm \hskip3cm \cdot\binom{2n-i-2k-j}{n-i-2k-j} a^{i+j}b^{n-i-2k-j}(c-b^2)^{k}  \nonumber\\
             \hskip-.22cm &=&\hskip-.22cm \sum_{i=0}^{n+1}(-1)^i\frac{i+1}{n+1}\sum_{k=0}^{n}\sum_{j=0}^{[\frac{n-k}{2}]}\binom{n+1}{k}\binom{n+1-k}{j}\\
             \hskip-.22cm & &\hskip-.22cm \hskip3cm  \cdot \binom{2n-i-k-2j}{n-i-k-2j} a^{i+k}b^{n-i-k-2j}(c-b^2)^{j}  \nonumber\\
             \hskip-.22cm &=&\hskip-.22cm \sum_{i=0}^{n+1}(-1)^i\frac{i+1}{n+1}\sum_{k=0}^{n}\sum_{j=0}^{n-k}\binom{n+1}{k}\binom{k}{j}\\
             \hskip-.22cm & &\hskip-.22cm \hskip3cm \cdot \binom{2n-i-k-j}{n-i-k-j} a^{i+k-j}b^{n-k-j}(c-b^2)^{j}.
\end{eqnarray*}
\end{proposition}

\section{A bijection between the sets $\mathcal{G}_n^{\mathbf{uvv}}(a,b,b^2)$ and $\mathcal{G}_n^{\mathbf{uvu}}(a,b,b^2)$  }

In this section, we give a direct bijection between the set $\mathcal{G}_n^{\mathbf{uvv}}(a,b,b^2)$ of $\mathbf{uvv}$-avoiding $(a,b,b^2)$-G-Motzkin paths and the set $\mathcal{G}_n^{\mathbf{uvu}}(a,b,b^2)$ of $\mathbf{uvu}$-avoiding $(a,b,b^2)$-G-Motzkin paths.

\begin{theorem}\label{theom 3.1.1}
There exists a bijection $\sigma$ between $\mathcal{G}_n^{\mathbf{uvv}}(a,b,b^2)$ and $\mathcal{G}_n^{\mathbf{uvu}}(a,b,b^2)$ for any integer $n\geq 0$.
\end{theorem}
\pf Given any $\mathbf{Q}\in \mathcal{G}_n^{\mathbf{uvv}}(a,b,b^2)$ for $n\geq 0$, when $n=0, 1$ and $2$, we define
$$\sigma(\varepsilon)=\varepsilon, \ \sigma(\mathbf{h}_a)=\mathbf{h}_a,\ \sigma(\mathbf{uv}_b)=\mathbf{uv}_b.$$

For $n\geq 2$, $\mathbf{Q}$ is $\mathbf{uvv}$-avoiding, there are six cases to be considered to define $\sigma(\mathbf{Q})$ recursively.

\subsection*{Case 1.} When $\mathbf{Q}=\mathbf{h}_a\mathbf{Q}'$ with $\mathbf{Q}'\in \mathcal{G}_{n-1}^{\mathbf{uvv}}(a,b,b^2)$, we define $\sigma(\mathbf{Q})=\mathbf{h}_a\sigma(\mathbf{Q}')$.

\subsection*{Case 2.} When $\mathbf{Q}=\mathbf{uv}_b\mathbf{h}_a\mathbf{Q}'$ with $\mathbf{Q}'\in \mathcal{G}_{n-2}^{\mathbf{uvv}}(a,b,b^2)$, we define $\sigma(\mathbf{Q})=\mathbf{uv}_b\mathbf{h}_a\sigma(\mathbf{Q}')$.

\subsection*{Case 3.} When $\mathbf{Q}=\mathbf{u}\mathbf{v}_b\mathbf{Q}''\mathbf{Q}'$ such that $\mathbf{Q}''\in \mathcal{G}_{k}^{\mathbf{uvv}}(a,b,b^2)$ is primitive and $\mathbf{Q}'\in \mathcal{G}_{n-k-1}^{\mathbf{uvv}}(a,b,b^2)$ for certain $1\leq k\leq n-1$, we define $\sigma(\mathbf{Q})=\mathbf{u}\sigma(\mathbf{Q}'')\mathbf{v}_b\sigma(\mathbf{Q}')$. In this case, one can notice that there exist $\mathbf{uvu}$'s in $\mathbf{Q}$, but not in $\sigma(\mathbf{Q})$.

\subsection*{Case 4.} When $\mathbf{Q}=\mathbf{u}^i\mathbf{ud}_{b^2}\mathbf{v}_b^i\mathbf{Q}'$ with $\mathbf{Q}'\in \mathcal{G}_{n-2-i}^{\mathbf{uvv}}(a,b,b^2)$ for $0\leq i\leq n-2$, we define
\begin{eqnarray*}
\sigma(\mathbf{Q})=\left\{
\begin{array}{rl}
 \mathbf{u}^{j}\mathbf{uv}_b\mathbf{d}_{b^2}^{j}\sigma(\mathbf{Q}'),  &  \mbox{if}\ i=2j-1\geq 1, \\[5pt]
 \mathbf{u}^{j+1}\mathbf{d}_{b^2}^{j+1}\sigma(\mathbf{Q}'),  & \mbox{if}\ i=2j\geq 0.
\end{array}\right.
\end{eqnarray*}

\subsection*{Case 5.} When $\mathbf{Q}=\mathbf{u}^i\mathbf{u}\mathbf{Q}''\mathbf{d}_{b^2}\mathbf{v}_b^i\mathbf{Q}'$ such that $\mathbf{Q}''\in \mathcal{G}_{k}^{\mathbf{uvv}}(a,b,b^2)$ is nonempty and $\mathbf{Q}'\in \mathcal{G}_{n-k-i}^{\mathbf{uvv}}(a,b,b^2)$ for certain $1\leq k\leq n-i$ and $0\leq i\leq n-2$, we define
\begin{eqnarray*}
\sigma(\mathbf{Q})=\left\{
\begin{array}{rl}
 \mathbf{u}^{j}\sigma(\mathbf{Q}''\mathbf{uv}_b)\mathbf{d}_{b^2}^{j}\sigma(\mathbf{Q}'),  &  \mbox{if}\ i=2j-1\geq 1, \\[5pt]
 \mathbf{u}^{j+1}\sigma(\mathbf{Q}''\mathbf{uv}_b)\mathbf{v}_b\mathbf{d}_{b^2}^{j}\sigma(\mathbf{Q}'),  & \mbox{if}\ i=2j\geq 0.
\end{array}\right.
\end{eqnarray*}

\subsection*{Case 6.} When $\mathbf{Q}=\mathbf{u}^i\mathbf{Q}''\mathbf{v}_b^i\mathbf{Q}'$ such that $\mathbf{Q}''\in \mathcal{G}_{k}^{\mathbf{uvv}}(a,b,b^2)$ is not primitive and $\mathbf{Q}'\in \mathcal{G}_{n-k-i}^{\mathbf{uvv}}(a,b,b^2)$ for certain $1\leq k\leq n-i$ and $1\leq i\leq n-1$, where $\mathbf{Q}''$ does not end with $\mathbf{uv}_b$ since $\mathbf{Q}$ is $\mathbf{uvv}$-avoiding, we define
\begin{eqnarray*}
\sigma(\mathbf{Q})=\left\{
\begin{array}{rl}
 \mathbf{u}^{j}\sigma(\mathbf{Q}'')\mathbf{v}_b\mathbf{d}_{b^2}^{j-1}\sigma(\mathbf{Q}'),  &  \mbox{if}\ i=2j-1\geq 1, \\[5pt]
 \mathbf{u}^{j}\sigma(\mathbf{Q}'')\mathbf{d}_{b^2}^{j}\sigma(\mathbf{Q}'),  & \mbox{if}\ i=2j\geq 2.
\end{array}\right.
\end{eqnarray*}

From the definition of $\sigma$, one can deduce by induction that $\sigma(\mathbf{Q})$ is $\mathbf{uvu}$-avoiding and the following assertations hold:

\begin{itemize}

\item In the case 3, $\sigma(\mathbf{Q}'')$ must be primitive and not be $\mathbf{uuv}_b\mathbf{v}_b$ since $\mathbf{Q}''$ is primitive;

\item In the case 5, $\sigma(\mathbf{Q}''\mathbf{uv}_b)$ has the form $\mathbf{P}_1\mathbf{uuv}_b\mathbf{v}_b$ or $\mathbf{P}_2\mathbf{uv}_b$ since $\mathbf{Q}''$ is nonempty, where both $\mathbf{P}_1\in \mathcal{G}_{r-1}^{\mathbf{uvu}}(a,b,b^2)$ and $\mathbf{P}_2\in \mathcal{G}_{r}^{\mathbf{uvu}}(a,b,b^2)$ must not end with $\mathbf{uv}_b$ for certain $r\geq 1$;

\item In the case 6, $\sigma(\mathbf{Q}'')$ is not primitive and does not end with $\mathbf{uv}_b$ or $\mathbf{uuv}_b\mathbf{v}_b$ since $\mathbf{Q}''$ is not primitive and does not end with $\mathbf{uv}_b$.

\end{itemize}

Conversely, the inverse procedure can be handled as follows. Given any $\mathbf{P}\in \mathcal{G}_{n}^{\mathbf{uvu}}(a,b,b^2)$ for $n\geq 0$, when $n=0, 1$, we define
$$\sigma^{-1}(\varepsilon)=\varepsilon, \ \sigma^{-1}(\mathbf{h}_a)=\mathbf{h}_a,\ \sigma^{-1}(\mathbf{uv}_b)=\mathbf{uv}_b.$$

For $n\geq 2$, there are five cases to be considered to define $\sigma^{-1}(\mathbf{P})$ recursively.

\subsection*{Case I} When $\mathbf{P}=\mathbf{h}_a\mathbf{P}'$ with $\mathbf{P}'\in \mathcal{G}_{n-1}^{\mathbf{uvu}}(a,b,b^2)$, we define $\sigma^{-1}(\mathbf{P})=\mathbf{h}_a\sigma^{-1}(\mathbf{P}')$.

\subsection*{Case II} When $\mathbf{P}=\mathbf{uv}_b\mathbf{h}_a\mathbf{P}'$ with $\mathbf{P}'\in \mathcal{G}_{n-2}^{\mathbf{uvu}}(a,b,b^2)$, we define $\sigma^{-1}(\mathbf{P})=\mathbf{uv}_b\mathbf{h}_a\sigma^{-1}(\mathbf{P}')$.

\subsection*{Case III} When $\mathbf{P}=\mathbf{u}\mathbf{P}''\mathbf{v}_b\mathbf{P}'$ such that $\mathbf{P}''\in \mathcal{G}_{k}^{\mathbf{uvu}}(a,b,b^2)$ and $\mathbf{P}'\in \mathcal{G}_{n-1-k}^{\mathbf{uvu}}(a,b,b^2)$ for certain $1\leq k\leq n-1$, we define
\begin{eqnarray*}
\sigma^{-1}(\mathbf{P})=\left\{
\begin{array}{rl}
 \mathbf{u}\mathbf{v}_b\sigma^{-1}(\mathbf{P}'')\sigma^{-1}(\mathbf{P}'),  &  \mbox{if $\mathbf{P}''$ is primitive and $\mathbf{P}''\neq \mathbf{uuv}_b\mathbf{v}_b $}, \\[5pt]
 \mathbf{u}\sigma^{-1}(\mathbf{P}'')\mathbf{v}_b\sigma^{-1}(\mathbf{P}'),  &  \mbox{if $\mathbf{P}''(\neq \varepsilon)$ is not primitive and does not} \\
                                                                           &  \ \ \ \mbox{ end with $\mathbf{uuv}_b\mathbf{v}_b$ or $\mathbf{uv}_b$}, \\[5pt]
 \mathbf{u}\sigma^{-1}(\mathbf{P}_1\mathbf{uv}_b)\mathbf{d}_{b^2}\sigma^{-1}(\mathbf{P}'),  &  \mbox{if $\mathbf{P}''=\mathbf{P}_1\mathbf{uuv}_b\mathbf{v}_b$}, \\[5pt]
 \mathbf{u}\sigma^{-1}(\mathbf{P}_2)\mathbf{d}_{b^2}\sigma^{-1}(\mathbf{P}'),  &  \mbox{if $\mathbf{P}''=\mathbf{P}_2\mathbf{uv}_b$},
\end{array}\right.
\end{eqnarray*}
where both $\mathbf{P}_1\in \mathcal{G}_{r-1}^{\mathbf{uvu}}(a,b,b^2)$ and $\mathbf{P}_2\in \mathcal{G}_{r}^{\mathbf{uvu}}(a,b,b^2)$ must not end with $\mathbf{uv}_b$ for certain $r\geq 1$, since $\mathbf{P}$ is $\mathbf{uvu}$-avoiding.

\subsection*{Case IV} When $\mathbf{P}=\mathbf{u}^{j}\mathbf{P}''\mathbf{d}_{b^2}^{j}\mathbf{P}'$ such that $\mathbf{P}''\in \mathcal{G}_{k}^{\mathbf{uvu}}(a,b,b^2)$ and $\mathbf{P}'\in \mathcal{G}_{n-2j-k}^{\mathbf{uvu}}(a,b,b^2)$ for certain $0\leq k\leq n-2j$ and the maximum $j\geq 1$, we define
\begin{eqnarray*}
\sigma^{-1}(\mathbf{P})=\left\{
\begin{array}{rl}
 \mathbf{u}^{2j-2}\mathbf{u}\mathbf{d}_{b^2}\mathbf{v}_{b}^{2j-2}\sigma^{-1}(\mathbf{P}'),  &  \mbox{if $\mathbf{P}''=\varepsilon$}, \\[5pt]
 \mathbf{u}^{2j-1}\mathbf{u}\sigma^{-1}(\mathbf{P}_1\mathbf{uv}_b)\mathbf{d}_{b^2}\mathbf{v}_{b}^{2j-1}\sigma^{-1}(\mathbf{P}'),  &  \mbox{if $\mathbf{P}''=\mathbf{P}_1\mathbf{uuv}_b\mathbf{v}_b$}, \\[5pt]
 \mathbf{u}^{2j-1}\mathbf{u}\sigma^{-1}(\mathbf{P}_2)\mathbf{d}_{b^2}\mathbf{v}_{b}^{2j-1}\sigma^{-1}(\mathbf{P}'),  &  \mbox{if $\mathbf{P}''=\mathbf{P}_2\mathbf{uv}_b$}, \\[5pt]
 \mathbf{u}^{2j}\sigma^{-1}(\mathbf{P}'')\mathbf{v}_b^{2j}\sigma^{-1}(\mathbf{P}'),  &  \mbox{if $\mathbf{P}''(\neq \varepsilon)$ is not primitive and does not} \\
                                                                           &  \ \ \ \mbox{ end with $\mathbf{uuv}_b\mathbf{v}_b$ or $\mathbf{uv}_b$},
\end{array}\right.
\end{eqnarray*}
where both $\mathbf{P}_1\in \mathcal{G}_{r-1}^{\mathbf{uvu}}(a,b,b^2)$ and $\mathbf{P}_2\in \mathcal{G}_{r}^{\mathbf{uvu}}(a,b,b^2)$ must not end with $\mathbf{uv}_b$ for certain $r\geq 1$, since $\mathbf{P}$ is $\mathbf{uvu}$-avoiding.

\subsection*{Case V} When $\mathbf{P}=\mathbf{u}^{j}\mathbf{u}\mathbf{P}''\mathbf{v}_b\mathbf{d}_{b^2}^{j}\mathbf{P}'$ such that $\mathbf{P}''\in \mathcal{G}_{k}^{\mathbf{uvu}}(a,b,b^2)$ and $\mathbf{P}'\in \mathcal{G}_{n-2j-1-k}^{\mathbf{uvu}}(a,b,b^2)$ for certain $0\leq k\leq n-2j-1$ and the maximum $j\geq 1$, we define
\begin{eqnarray*}
\sigma^{-1}(\mathbf{P})=\left\{
\begin{array}{rl}
 \mathbf{u}^{2j-1}\mathbf{u}\mathbf{d}_{b^2}\mathbf{v}_{b}^{2j-1}\sigma^{-1}(\mathbf{P}'),  &  \mbox{if $\mathbf{P}''=\varepsilon$}, \\[5pt]
 \mathbf{u}^{2j}\mathbf{u}\sigma^{-1}(\mathbf{P}_1\mathbf{uv}_b)\mathbf{d}_{b^2}\mathbf{v}_{b}^{2j}\sigma^{-1}(\mathbf{P}'),  &  \mbox{if $\mathbf{P}''=\mathbf{P}_1\mathbf{uuv}_b\mathbf{v}_b$}, \\[5pt]
 \mathbf{u}^{2j}\mathbf{u}\sigma^{-1}(\mathbf{P}_2)\mathbf{d}_{b^2}\mathbf{v}_{b}^{2j}\sigma^{-1}(\mathbf{P}'),  &  \mbox{if $\mathbf{P}''=\mathbf{P}_2\mathbf{uv}_b$}, \\[5pt]
 \mathbf{u}^{2j+1}\sigma^{-1}(\mathbf{P}'')\mathbf{v}_b^{2j+1}\sigma^{-1}(\mathbf{P}'),  &  \mbox{if $\mathbf{P}''(\neq \varepsilon)$ is not primitive and does not} \\
                                                                           &  \ \ \ \mbox{ end with $\mathbf{uuv}_b\mathbf{v}_b$ or $\mathbf{uv}_b$},
\end{array}\right.
\end{eqnarray*}
where both $\mathbf{P}_1\in \mathcal{G}_{r-1}^{\mathbf{uvu}}(a,b,b^2)$ and $\mathbf{P}_2\in \mathcal{G}_{r}^{\mathbf{uvu}}(a,b,b^2)$ must not end with $\mathbf{uv}_b$ for certain $r\geq 1$, since $\mathbf{P}$ is $\mathbf{uvu}$-avoiding.

It is not difficult to verify that $\sigma^{-1}\sigma=\sigma\sigma^{-1}=1$, both $\sigma$ and $\sigma^{-1}$ are two weight-keeping mappings and $\sigma^{-1}(\mathbf{P})$ is
$\mathbf{uvv}$-avoiding by induction on the length of $\mathbf{P}$. Hence, $\sigma$ is a desired bijection between $\mathcal{G}_n^{\mathbf{uvv}}(a,b,b^2)$ and $\mathcal{G}_n^{\mathbf{uvu}}(a,b,b^2)$. This completes the proof of Theorem \ref{theom 3.1.1}. \qed\vskip0.2cm

In order to give a more intuitive view on the bijection $\sigma$, a pictorial description of $\sigma$ is presented for $\mathbf{Q}=\mathbf{u}^3\mathbf{d}_{b^2}\mathbf{v}_b^2\mathbf{u}^2\mathbf{d}_{b^2}\mathbf{v}_b
\mathbf{u}^5\mathbf{v}_b\mathbf{d}_{b^2}\mathbf{v}_b^3\mathbf{h}_{a}
\mathbf{u}^2\mathbf{h}_{a}\mathbf{d}_{b^2}\mathbf{v}_b\mathbf{u}^3\mathbf{v}_b\mathbf{h}_a\mathbf{v}_b\mathbf{u}\mathbf{v}_b
\mathbf{u}^2\mathbf{h}_a\mathbf{u}\mathbf{d}_{b^2}\mathbf{v}_b^3 \in \mathcal{G}_{28}^{\mathbf{uvv}}(a,b,b^2)$, we have
$$\small\sigma(\mathbf{Q})=\mathbf{u}^2\mathbf{d}_{b^2}^2\mathbf{u}^2\mathbf{v}_b\mathbf{d}_{b^2}
\mathbf{u}^4\mathbf{v}_b^2\mathbf{d}_{b^2}^2\mathbf{h}_{a}
\mathbf{u}\mathbf{h}_{a}\mathbf{u}\mathbf{v}_b\mathbf{d}_{b^2}\mathbf{u}^3\mathbf{v}_b\mathbf{h}_a\mathbf{v}_{b}\mathbf{u}^2\mathbf{h}_{a} \mathbf{u}\mathbf{d}_{b^2}^2\mathbf{v}_b^2\in \mathcal{G}_{28}^{\mathbf{uvu}}(a,b,b^2). $$

See Figure 2 for detailed illustrations.

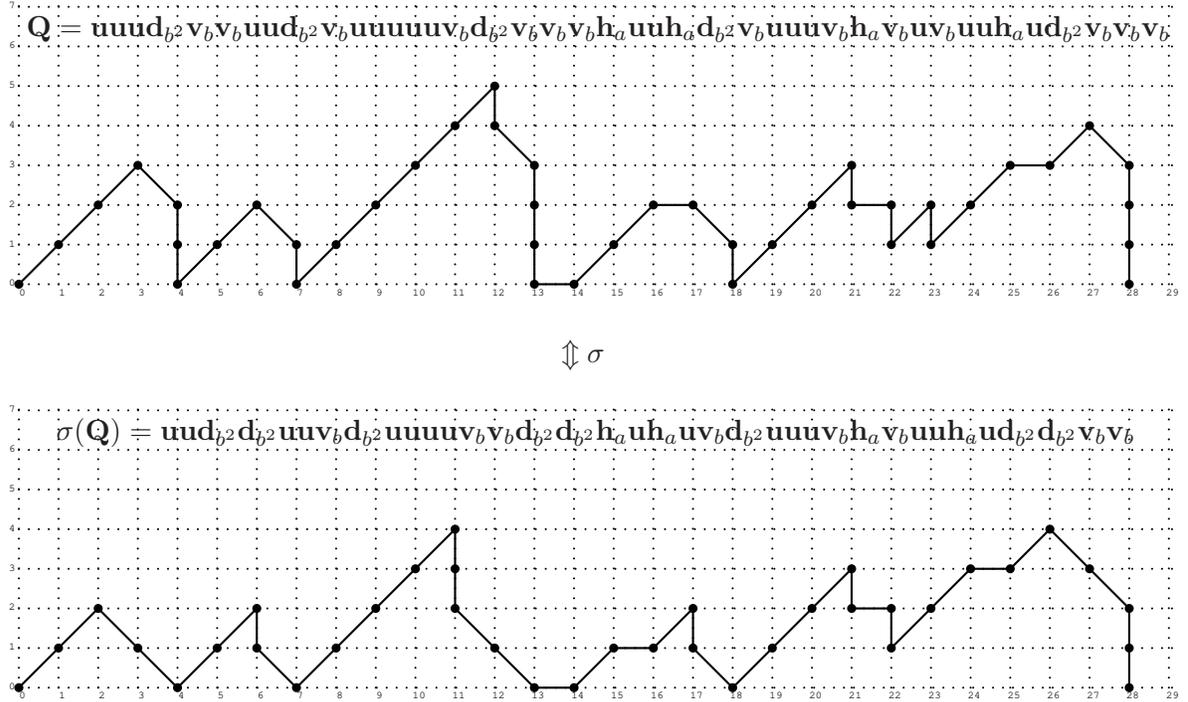
\begin{figure}[h] \setlength{\unitlength}{0.5mm}

\begin{center}
\begin{pspicture}(15,4)
\psset{xunit=15pt,yunit=15pt}\psgrid[subgriddiv=1,griddots=4,
gridlabels=4pt](0,0)(29,7)

\psline(0,0)(3,3)(4,2)(4,0)(6,2)(7,1)(7,0)(12,5)(12,4)(13,3)(13,0)(14,0)(16,2)(17,2)(18,1)(18,0)(21,3)(21,2)(22,2)(22,1)(23,2)(23,1)(25,3)(26,3)(27,4)(28,3)(28,0)

\pscircle*(0,0){0.06}\pscircle*(1,1){0.06}\pscircle*(2,2){0.06}
\pscircle*(3,3){0.06}\pscircle*(4,0){0.06}\pscircle*(4,1){0.06}\pscircle*(4,2){0.06}\pscircle*(5,1){0.06}
\pscircle*(6,2){0.06}\pscircle*(7,1){0.06}\pscircle*(7,0){0.06}

\pscircle*(8,1){0.06}\pscircle*(9,2){0.06}\pscircle*(10,3){0.06}\pscircle*(11,4){0.06}
\pscircle*(12,5){0.06}\pscircle*(12,4){0.06}\pscircle*(13,3){0.06}
\pscircle*(13,2){0.06}\pscircle*(13,1){0.06}\pscircle*(13,0){0.06}\pscircle*(14,0){0.06}
\pscircle*(15,1){0.06}\pscircle*(16,2){0.06}\pscircle*(17,2){0.06}\pscircle*(18,1){0.06}\pscircle*(18,0){0.06}
\pscircle*(19,1){0.06}\pscircle*(20,2){0.06}\pscircle*(21,3){0.06}\pscircle*(21,2){0.06}\pscircle*(22,2){0.06}\pscircle*(22,1){0.06}
\pscircle*(23,2){0.06}\pscircle*(23,1){0.06}\pscircle*(24,2){0.06}\pscircle*(25,3){0.06}\pscircle*(26,3){0.06}\pscircle*(27,4){0.06}
\pscircle*(28,3){0.06}\pscircle*(28,2){0.06}\pscircle*(28,1){0.06}\pscircle*(28,0){0.06}

\put(.1,3.3){$\mathbf{Q}=\mathbf{uuu}\mathbf{d}_{b^2}\mathbf{v}_b\mathbf{v}_b\mathbf{uu}\mathbf{d}_{b^2}\mathbf{v}_b
\mathbf{uuuuu}\mathbf{v}_b\mathbf{d}_{b^2}\mathbf{v}_b\mathbf{v}_b\mathbf{v}_b\mathbf{h}_{a} 
\mathbf{uu}\mathbf{h}_{a}\mathbf{d}_{b^2}\mathbf{v}_b\mathbf{uuu}\mathbf{v}_b\mathbf{h}_a\mathbf{v}_b\mathbf{u}\mathbf{v}_b
\mathbf{uu}\mathbf{h}_a\mathbf{u}\mathbf{d}_{b^2}\mathbf{v}_b\mathbf{v}_b\mathbf{v}_b$}

\end{pspicture}
\end{center}\vskip0.5cm

$\Updownarrow \sigma$
\vskip0.2cm

\begin{center}
\begin{pspicture}(15,4)
\psset{xunit=15pt,yunit=15pt}\psgrid[subgriddiv=1,griddots=4,
gridlabels=4pt](0,0)(29,7)

\psline(0,0)(2,2)(4,0)(6,2)(6,1)(7,0)(11,4)(11,2)(13,0)(14,0)(15,1)(16,1)(17,2)(17,1)(18,0)(21,3)(21,2)(22,2)(22,1)(24,3)(25,3)(26,4)(28,2)(28,0)

\pscircle*(0,0){0.06}\pscircle*(1,1){0.06}\pscircle*(2,2){0.06}
\pscircle*(3,1){0.06}\pscircle*(4,0){0.06}\pscircle*(5,1){0.06}
\pscircle*(6,2){0.06}\pscircle*(6,1){0.06}\pscircle*(7,0){0.06}

\pscircle*(8,1){0.06}\pscircle*(9,2){0.06}\pscircle*(10,3){0.06}\pscircle*(11,4){0.06}\pscircle*(11,3){0.06}\pscircle*(11,2){0.06}
\pscircle*(12,1){0.06}\pscircle*(13,0){0.06}\pscircle*(14,0){0.06}\pscircle*(15,1){0.06}
\pscircle*(16,1){0.06}\pscircle*(17,2){0.06}\pscircle*(17,1){0.06}\pscircle*(18,0){0.06}
\pscircle*(19,1){0.06}\pscircle*(20,2){0.06}\pscircle*(21,3){0.06}\pscircle*(21,2){0.06}
\pscircle*(22,2){0.06}\pscircle*(22,1){0.06}\pscircle*(23,2){0.06}\pscircle*(24,3){0.06}\pscircle*(25,3){0.06}\pscircle*(26,4){0.06}
\pscircle*(27,3){0.06}\pscircle*(28,2){0.06}\pscircle*(28,1){0.06}\pscircle*(28,0){0.06}

\put(.5,3.3){$\sigma(\mathbf{Q})=\mathbf{uu}\mathbf{d}_{b^2}\mathbf{d}_{b^2}\mathbf{uu}\mathbf{v}_b\mathbf{d}_{b^2}
\mathbf{uuuu}\mathbf{v}_b\mathbf{v}_b\mathbf{d}_{b^2}\mathbf{d}_{b^2}\mathbf{h}_{a}
\mathbf{u}\mathbf{h}_{a}\mathbf{u}\mathbf{v}_b\mathbf{d}_{b^2}\mathbf{uuu}\mathbf{v}_b\mathbf{h}_a\mathbf{v}_{b}\mathbf{uu}\mathbf{h}_{a} \mathbf{u}\mathbf{d}_{b^2}\mathbf{d}_{b^2}\mathbf{v}_b\mathbf{v}_b$}

\end{pspicture}
\end{center}\vskip0.3cm

\caption{\small An example of the bijection $\sigma$ described in the proof of Theorem \ref{theom 3.1.1}. }

\end{figure}

\section{Counting the set of fixed points of the bijection $\sigma$ }

In this section, we will count the set of fixed points of the bijection $\sigma$ presented in Section 3.

Let $\mathcal{F}_n=\{\mathbf{Q}\in \mathcal{G}_{n}^{\mathbf{uvv}}(a,b,b^2)|\sigma(\mathbf{Q})=\mathbf{Q}\}$ and $\mathcal{F}=\bigcup_{n\geq 0}\mathcal{F}_n$, set $F_n=|\mathcal{F}_n|$. It is easy to verify the few initial values for $F_n$, see Table 4.1.
\begin{center}
\begin{eqnarray*}
\begin{array}{c|ccccccccccc}\hline
   n    & 0   & 1   &  2   & 3    & 4    & 5    &  6     &   7   &    8   &     9     &     10      \\\hline
  F_n   & 1   & 2   &  5   & 13   & 39   & 125  &  421   &  1478 &  5329  &   19658   & 73783       \\\hline
\end{array}
\end{eqnarray*}
Table 4.1. The first values of $F_n$.
\end{center}

According to the definition of $\sigma$, any $\mathbf{Q}\in \mathcal{F}_n$ must belong to $\mathcal{G}_{n}^{\{\mathbf{uvv}, \mathbf{uvu}\}}(a,b,b^2)$, the set of $(a, b, b^2)$-G-Motzkin paths avoiding both the $\mathbf{uvv}$ and $\mathbf{uvu}$ patterns, since $\mathbf{Q}$ is $\mathbf{uvv}$-avoiding and $\sigma(\mathbf{Q})$ is $\mathbf{uvu}$-avoiding. But there exists $\mathbf{P}\in \mathcal{G}_{n}^{\{\mathbf{uvv}, \mathbf{uvu}\}}(a,b,b^2)$ such that $\sigma(\mathbf{P})\neq \mathbf{P}$. For example, $\sigma(\mathbf{uudv})=\mathbf{uuvd}\neq \mathbf{uudv}$ for $n=3$. Precisely, in the proof of Theorem \ref{theom 3.1.1}, one can deduce that $\mathbf{Q}\notin \mathcal{F}_n(\sigma)$ if $\mathbf{Q}$ is being in the following situations, 1) in the whole case 3; 2) in the case 4 when $i\geq 1$; 3) in the whole case 5; and 4) in the case 6 when $i\geq 2$.
Equivalently, one can derive that
\begin{itemize}

\item In Case 1, $\mathbf{Q}=\mathbf{h}_a\mathbf{Q}'\in \mathcal{F}_n$ if and only if $\mathbf{Q}'\in \mathcal{F}_{n-1}$;\vskip0.1cm

\item In Case 2, $\mathbf{Q}=\mathbf{uv}_b\mathbf{h}_a\mathbf{Q}'\in \mathcal{F}_n$ if and only if $\mathbf{Q}'\in \mathcal{F}_{n-2}$;\vskip0.1cm

\item In Case 4 when $i=0$, $\mathbf{Q}=\mathbf{ud}_{b^2}\mathbf{Q}'\in \mathcal{F}_n$ if and only if $\mathbf{Q}'\in \mathcal{F}_{n-2}$;\vskip0.1cm

\item In Case 6 when $i=1$, $\mathbf{Q}=\mathbf{u}\mathbf{Q}''\mathbf{v}_{b}\mathbf{Q}'\in \mathcal{F}_n$ if and only if $\mathbf{Q}''\in \mathcal{F}_{k}$ and $\mathbf{Q}'\in \mathcal{F}_{n-k-1}$ for certain $1\leq k\leq n-1$ such that $\mathbf{Q}''(\neq \varepsilon)$ is not primitive and does not end with $\mathbf{uv}_b$.

\end{itemize}

Let $\mathcal{A}_n$ be the subset of $\mathbf{Q}\in \mathcal{F}_n$ such that $\mathbf{Q}$ is not primitive and does not end with $\mathbf{uv}_b$, $\mathcal{B}_n$ be the subset of $\mathbf{Q}\in \mathcal{F}_n$ such that $\mathbf{Q}$ ends with $\mathbf{uv}_b$, and $\mathcal{C}_n$ be the subset of $\mathbf{Q}\in \mathcal{F}_n$ such that $\mathbf{Q}$ is primitive and does not end with $\mathbf{uv}_b$. Set $a_n=|\mathcal{A}_n|, b_n=|\mathcal{B}_n|, c_n=|\mathcal{C}_n|$. Firstly, $\mathcal{F}_n$ is the disjoint union of $\mathcal{A}_n, \mathcal{B}_n$ and $\mathcal{C}_n$, i.e., $F_n=a_n+b_n+c_n$ for $n\geq 0$; Secondly, $\mathcal{C}_0=\mathcal{C}_1=\emptyset$, $\mathcal{C}_2=\{\mathbf{uh}_a\mathbf{v}_b, \mathbf{ud}_{b^2}\}$ and $\mathcal{C}_n=\mathbf{u}\mathcal{A}_{n-1}\mathbf{v}_b$ for $n\geq 3$, i.e., $c_n=a_{n-1}$ for $n\geq 3$ with $c_0=c_1=0$ and $c_2=c_3=2$; Thirdly, $\mathcal{B}_n$ is the disjoint union of $\mathcal{A}_{n-1}\mathbf{uv}_b$ and $\mathcal{C}_{n-1}\mathbf{uv}_b$ for $n\geq 1$, i.e., $b_n=a_{n-1}+c_{n-1}$ for $n\geq 1$ with $b_0=0$ and $b_1=b_2=1$.
These together generate the following Lemma.

\begin{lemma} For any integer $n\geq 4$, there holds
\begin{eqnarray}\label{eqn 4.1}
F_n \hskip-.22cm &=&\hskip-.22cm a_n+2a_{n-1}+a_{n-2}
\end{eqnarray}
with $a_0=a_1=1, a_2=2, a_3=7$ and $a_4=23$.
\end{lemma} \vskip0.1cm

On the other hand, the family $\mathcal{F}$ can be partitioned into the form:
$$\mathcal{F}=\varepsilon+\mathbf{h}_a\mathcal{F}+\mathbf{uv}_b\mathbf{h}_a\mathcal{F}+ \mathbf{ud}_{b^2}\mathcal{F}+\mathbf{u}\mathcal{A}'\mathbf{v}_b\mathcal{F}, $$
where $\mathcal{A}'=\mathcal{A}-\varepsilon$ and $\mathcal{A}=\bigcup_{n\geq 0}\mathcal{A}_n$. This leads to the following recurrence for $F_n$.

\begin{lemma} For any integer $n\geq 1$, there holds
\begin{eqnarray}\label{eqn 4.2}
F_{n+1} \hskip-.22cm &=&\hskip-.22cm F_n+2F_{n-1}+\sum_{k=1}^{n}a_kF_{n-k}
\end{eqnarray}
with $F_0=1, F_1=2$.
\end{lemma} \vskip0.1cm

Let $F(x)=\sum_{n\geq 0}F_nx^n$ and $A(x)=\sum_{n\geq 0}a_nx^n$. By (\ref{eqn 4.1}), we have
\begin{eqnarray}\label{eqn 4.3}
F(x) \hskip-.22cm &=&\hskip-.22cm 1+2x+5x^2+13x^3+\sum_{n\geq 4}(a_n+2a_{n-1}+a_{n-2})x^n \nonumber \\
     \hskip-.22cm &=&\hskip-.22cm 1+2x+5x^2+13x^3+\big(A(x)-1-x-2x^2-7x^3\big)            \nonumber \\
     \hskip-.22cm & &\hskip-.22cm \ \ +\ 2x\big(A(x)-1-x-2x^2\big)+x^2\big(A(x)-1-x\big)  \nonumber \\
     \hskip-.22cm &=&\hskip-.22cm (1+x)^2A(x)-x+x^3.
\end{eqnarray}
By (\ref{eqn 4.2}), we can obtain
\begin{eqnarray}\label{eqn 4.4}
F(x) \hskip-.22cm &=&\hskip-.22cm 1+2x+x(F(x)-1)+2x^2F(x)+x(A(x)-1)F(x)  \nonumber \\
     \hskip-.22cm &=&\hskip-.22cm 1+x+2x^2F(x)+xA(x)F(x).
\end{eqnarray}

Eliminating $A(x)$ in (\ref{eqn 4.3}) and (\ref{eqn 4.4}) produces
\begin{eqnarray*}
xF(x)^2-(1+x)(1+x-3x^2-x^3)F(x)+(1+x)^3 \hskip-.22cm &=&\hskip-.22cm 0.
\end{eqnarray*}
Solve this, we have
\begin{eqnarray}\label{eqn 4.5}
F(x) \hskip-.22cm &=&\hskip-.22cm \frac{(1+x)(1+x-3x^2-x^3)-(1+x)\sqrt{(1+x-3x^2-x^3)^2-4x(1+x)}}{2x} \nonumber \\
     \hskip-.22cm &=&\hskip-.22cm  \frac{(1+x)^2}{1+x-3x^2-x^3}C\Big(\frac{x(1+x)}{(1+x-3x^2-x^3)^2}\Big).
\end{eqnarray}
By (\ref{eqn 4.5}), taking the coefficient of $x^n$ in $F(x)$, we get the explicit formula for the number $F_n$ of the fixed points of the bijection $\sigma$, namely,
\begin{theorem}
For any integer $n\geq 0$, there holds
\begin{eqnarray*}
F_n \hskip-.22cm &=&\hskip-.22cm \sum_{k=0}^{n}\sum_{j=0}^{[\frac{n-k}{2}]}\sum_{i=0}^{j}(-1)^{n-k-i}\binom{2k+j}{j}\binom{j}{i}\binom{n-j-i-2}{n-k-2j-i}3^{j-i}C_k,
\end{eqnarray*}
where $C_k$ is the $k$-th Catalan number.
\end{theorem}

\vskip0.5cm
\section*{Declaration of competing interest}

The authors declare that they have no known competing financial interests or personal relationships that could have
appeared to influence the work reported in this paper.

\section*{Acknowledgements} {The authors are grateful to the referees for
the helpful suggestions and comments. The Project is sponsored by ``Liaoning
BaiQianWan Talents Program". }

\vskip.2cm


\end{document}